\documentclass[12pt]{article}
\usepackage{latexsym,amsmath}
\usepackage{url}
\usepackage{amssymb}

\newtheorem{thm}{Theorem}[section]

\newtheorem{lemma}[thm]{Lemma}
\newtheorem{prop}[thm]{Proposition}
\newtheorem{cor}[thm]{Corollary}

\newtheorem{remark}[thm]{Remark}

\newcommand{\comp}{{\operatorname{\rm comp}}}

\newcommand{\beq}[1]{\begin{equation}\label{#1}}
\newcommand{\enq}[0]{\end{equation}}

\newcommand{\qed}[0]{{\hspace*{\fill}\mbox{$\Box$}}}

\newcommand{\cA}[0]{{\cal A}}
\newcommand{\cC}[0]{{\cal C}}
\newcommand{\cE}[0]{{\cal E}}
\newcommand{\cF}[0]{{\cal F}}
\newcommand{\cH}[0]{{\cal H}}

\newcommand{\cM}[0]{{\cal M}}
\newcommand{\cO}[0]{{\cal O}}
\newcommand{\cR}[0]{{\cal R}}

\newcommand{\cU}[0]{{\cal U}}

\newcommand{\cW}[0]{{\cal W}}
\newcommand{\cX}[0]{{\cal X}}

\newcommand{\Z}{{\mathbb Z}}

\newcommand{\ga}[0]{\alpha}
\newcommand{\gb}[0]{\beta}
\newcommand{\gre}[0]{\varepsilon}
\newcommand{\gd}[0]{\delta}
\newcommand{\gD}[0]{\Delta}

\newcommand{\gl}[0]{\lambda}
\newcommand{\go}[0]{\omega}
\newcommand{\gO}[0]{\Omega}
\newcommand{\gS}[0]{\Sigma}

\begin{document}

\renewcommand{\thefootnote}{\fnsymbol{footnote}}
\footnotetext{Key words: Mixing time, 3-colouring, Potts model,
conductance, Glauber dynamics, discrete hypercube.}
\footnotetext{Mathematics Subject Classifications: 05C15, 82B20.}

\title{Sampling $3$-colourings of regular bipartite graphs}

\author{David Galvin\thanks{Department of Mathematics,
University of Pennsylvania, 209 South 33rd Street, Philadelphia PA
19104; dgalvin@math.upenn.edu. This work was begun while the author
was a member of the Institute for Advanced Study, Einstein Drive,
Princeton, NJ 08540 and was supported in part by NSF grant
DMS-0111298.}}

\date{Submitted August 7, 2006; accepted April 5, 2007}

\maketitle

\begin{abstract}

We show that if $\gS=(V,E)$ is a regular bipartite graph for which
the expansion of subsets of a single parity of $V$ is reasonably
good and which satisfies a certain local condition (that the union
of the neighbourhoods of adjacent vertices does not contain too many
pairwise non-adjacent vertices), and if $\cM$ is a Markov chain on
the set of proper $3$-colourings of $\gS$ which updates the colour
of at most $\rho|V|$ vertices at each step and whose stationary
distribution is uniform, then for $\rho \approx .22$ and $d$
sufficiently large the convergence to stationarity of $\cM$ is
(essentially) exponential in $|V|$. In particular, if $\gS$ is the
$d$-dimensional hypercube $Q_d$ (the graph on vertex set $\{0,1\}^d$
in which two strings are adjacent if they differ on exactly one
coordinate) then the convergence to stationarity of the well-known
Glauber (single-site update) dynamics is exponentially slow in
$2^d/(\sqrt{d}\log d)$. A combinatorial corollary of our main result
is that in a uniform $3$-colouring of $Q_d$ there is an
exponentially small probability (in $2^d$) that there is a colour
$i$ such the proportion of vertices of the even subcube coloured $i$
differs from the proportion of the odd subcube coloured $i$ by at
most $.22$. Our proof combines a conductance argument with
combinatorial enumeration methods.
\end{abstract}

\newpage

\section{Introduction and statement of the result}
\label{sec-intro}

Markov chain Monte Carlo algorithms (MCMC's) occur frequently in
computer science in algorithms designed to sample from or estimate
the size of large combinatorially defined structures; they are also
used in statistical physics and the study of networks to help
understand the behavior of models of physical systems and networks
in equilibrium. In this paper we study a class of natural MCMC's
that sample from proper $3$-colourings of a regular bipartite graph.

Let $\gS=(V,E)$ be a simple, loopless, finite graph on vertex set
$V$ and edge set $E$. (For graph theory basics, see {\em e.g.}
\cite{Bollobas}, \cite{Diestel}.) For a positive integer $q$ write
$\cC_q=\cC_q(\gS)$ for the set of proper $q$-colourings of $\gS$;
that is,
$$
\cC_q=\{\chi:V(\gS)\rightarrow \{0,1,\ldots,q-1\}:xy \in
E(\gS)\Rightarrow\chi(x)\neq\chi(y)\}.
$$
Let $\pi_q=\pi_q(\gS)$ be the uniform probability distribution on
$\cC_q$.

The notion of $q$-colouring is fundamental in graph theory; see {\em
e.g.} \cite[Chapter 5]{Bollobas2} for a survey. The notion also
occurs in statistical physics; the pair $(\cC_q,\pi_q)$ is the
zero-temperature limit of the $q$-state antiferromagnetic Potts
model (see {\em e.g.} \cite{Sokal1,Sokal2}).

\medskip

{\em Glauber dynamics} for proper $q$-colourings is the single-site
update Markov chain $\cM_q=\cM_q(\gS)$ on state space $\cC_q$ with
transition probabilities $P_q(\chi_1,\chi_2)$, $\chi_1,\chi_2 \in
\cC_q,$ given by
$$
P_q(\chi_1,\chi_2) = \left\{
            \begin{array}{ll}
               0 & \mbox{ if $|\{v \in V:\chi_1(v)\neq \chi_2(v)\}| > 1$} \\
               \frac{1}{|V|}\frac{1}{q} & \mbox{ if $|\{v \in V:\chi_1(v)\neq \chi_2(v)\}| = 1$} \\
               1 - \sum_{\chi_1 \neq \chi_2' \in \cC_q} P_q(\chi_1,\chi_2') & \mbox{ if
$\chi_1=\chi_2$.}
            \end{array}
         \right.
$$
We may think of $\cM_q$ dynamically as follows. From a $q$ colouring
$\chi$, choose a vertex $v$ uniformly from $V$ and a colour $j$
uniformly from $\{0,\ldots,q-1\}$. Then define a function
$\chi':V\rightarrow \{0,\ldots,q-1\}$ by
$$
\chi'(w)=\left\{ \begin{array}{ll}
                    \chi(w) & ~\mbox{if $w \neq v$} \\
                    j & ~\mbox{if $w = v$.}
                 \end{array}
          \right.
$$
Finally, move to $\chi'$ if $\chi'$ is a proper $q$-colouring, and
stay at $\chi$ otherwise. (A variant of Glauber dynamics chooses $j$
uniformly from $\{0,\ldots,q-1\}\setminus \{\chi(w):w \sim v\}$,
ensuring that $\chi'$ is always a proper colouring. This changes the
transition probabilities, but does not significantly change the
qualitative behavior of the chain.)

For all $\gS$ the chain $\cM_q$ is aperiodic, but it is not in
general irreducible (consider, for example, $\gS=K_q$, the complete
graph on $q$ vertices), and so not ergodic. In the case when $\cM_q$
is ergodic ({\em e.g.}, when $\gS$ has maximum degree $\gD$ and $q
\geq \gD+2$; see \cite{Jerrum}) it is readily checked that it has
(unique) stationary distribution $\pi_q$. (One only has to check
that $\cM_q$ is reversible with respect to $\pi_q$; that is, that it
satisfies the detailed balance equations
$\pi_q(\chi_1)P_q(\chi_1,\chi_2)=\pi_q(\chi_2)P_q(\chi_2,\chi_1)$
for all $\chi_1, \chi_2 \in \cC_q$.) A natural and important
question to ask about $\cM_q$ in this case is how quickly it
converges to its stationary distribution. We define the {\em mixing
time} $\tau_{{\cal M}_q}$ of ${\cal M}_q$ by
$$
\tau_{\cM_q}=\min \left\{t~:~ d_{TV}(P_q^t,\pi_q) \leq
\frac{1}{e}\right\}
$$
where $P_q^t(\chi,\chi')$ is the probability of moving from $\chi$
to $\chi'$ in $t$ steps and
$$
d_{TV}(P_q^t,\pi_q) = \max_{\chi_1 \in \cC_q} \frac{1}{2}
\sum_{\chi_2 \in \cC_q} |P_q^t(\chi_1,\chi_2)-\pi_q(\chi_2)|
$$
is total variation distance. The mixing time of $\cM_q$ captures the
speed at which the chain converges to its stationary distribution:
for every $\gre
>0$, in order to get a sample from $\cC_q$ which is within
$\gre$ of $\pi_q$ (in total variation distance), it is necessary and
sufficient to run the chain from some arbitrarily chosen
distribution for some multiple (depending on $\gre$) of the mixing
time. For surveys of issues related to the mixing time of a Markov
chain, see {\em e.g.} \cite{AldousFill,MontenegroTetali,Randall}.

Jerrum \cite{Jerrum} and Salas and Sokal \cite{SalasSokal}
independently showed that if $\gS$ has maximum degree $\gD$ and
$q>2\gD$ then there is rapid mixing of the Glauber dynamics; {\em
i.e.}, $\tau_{\cM_q}(\gS)$ is polynomial in $|V|$. In fact, they
showed that the mixing time is optimal ($O(|V|\log |V|)$). Bubley
and Dyer \cite{BubleyDyer} showed that there is rapid mixing for
$q=2\gD$ and Molloy \cite{Molloy} improved this to optimal mixing.
In a breakthrough result Vigoda \cite{Vigoda} showed rapid mixing
for $q \geq (11/6)\gD$. More recently Dyer, Greenhill and Molloy
\cite{DyerGreenhillMolloy} exhibited optimal mixing for $q \geq
(2-\gre)\Delta$ for a small positive constant $\gre$.

\medskip

In this paper our aim is to explore the limitations of Glauber
dynamics as a sampling tool by exhibiting a class of graphs for
which the mixing time is essentially as far from optimal as
possible. In this direction, \L uczak and Vigoda \cite{LuczakVigoda}
have exhibited families of planar graphs for which $\cM_q$ is not
rapidly mixing for each fixed $q\geq 3$ and families of bipartite
graphs with maximum degree $\gD$ for which $\cM_q$ is not rapidly
mixing for any $3 \leq q \leq O(\gD/\log \gD)$. A drawback of these
negative results is that the families exhibited consist of random
graphs. Here, we attempt to remedy this by constructing explicit
families of graphs for which Glauber dynamics is inefficient. We
focus exclusively on the case $q=3$ (we cannot see at the moment how
to apply our techniques to any $q>3$) and $\gS$ regular bipartite.
Specifically, we establish certain local and expansion conditions in
a regular bipartite graph $\gS$ that force $\tau_{\cM_3(\gS)}$ to be
(almost) exponential in $|V|$. The discrete hypercube is among the
families of graphs which satisfy our conditions.

Our techniques actually apply to the class of $\rho$-local chains
(considered in \cite{BorgsChayesFriezeKimTetaliVigodaVu} and also in
\cite{DyerFriezeJerrum}, where the terminology {\em
$\rho|V|$-cautious} is employed) for suitably small $\rho$. A Markov
chain $\cM$ on state space $\cC_q$ is {\em $\rho$-local} if in each
step of the chain at most $\rho|V|$ vertices have their colour
changed; that is, if
$$
P_\cM(\chi_1,\chi_2) \neq 0 \Rightarrow |\{v \in V : \chi_1(v)\neq
\chi_2(v)\}| \leq \rho |V|.
$$

Before stating our main result, we establish some notation. From now
on, $\gS=(V,E)$ will be a $d$-regular bipartite graph with partition
classes $\cE$ and $\cO$. For $u, v \in V$ we write $u \sim v$ if
there is an edge in $\gS$ joining $u$ and $v$. Set $N(u)=\{w\in V:w
\sim u\}$ ($N(u)$ is the {\em neighbourhood} of $u$) and for $A
\subseteq V$ set $N(A)=\cup_{w \in A} N(w)$. For $A \subseteq \cE$
(or $A \subseteq \cO$) set
$$
[A]=\{x \in V : N(x) \subseteq N(A)\}
$$
(we think of $[A]$ as an external closure of $A$) and say that such
an $A$ is {\em small} if $|[A]| \leq |V|/4$. Note that $N(A)$
determines $[A]$ but not $A$ itself.

Define the {\em bipartite expansion} of $\gS$ by
$$
\gd(\gS)=\min \left\{\frac{|N(A)|-|[A]|}{|N(A)|}:A \subseteq \cE
~\mbox{(or $\cO$)}~\mbox{small},~A \neq \emptyset\right\};
$$
note that $0 \leq \gd < 1$. The second inequality is clear. To see
the first, note that since $\gS$ is regular and bipartite it has a
perfect matching, and so satisfies
\begin{equation} \label{obs-bipartite,regular}
\mbox{$|X| \leq |N(X)|$ for all $X \subseteq \cE$ or $\cO$.}
\end{equation}
That $0 \leq \gd$ now follows from $|[A]|\leq |N([A])|=|N(A)|$. The
bipartite expansion constant is a measure of the proportion by which
the neighbourhood size of a small set exceeds the size of the set
itself, in the worst case.

Finally, define the {\em locality} $\ell(\gS)$ of $\gS$ to be the
largest $\ell \geq 0$ such that for all $x \sim y \in V$ and for all
independent sets $I$ (sets of vertices spanning no edges) in the
subgraph of $\gS$ induced by $N(x) \cup N(y)$ we have $|I| \leq
2d-\ell$. (So, for example, if $\gS$ is the $d$-regular tree then
$\ell(\gS)=2$ since the subgraph induced by the neighbourhoods of
adjacent vertices contains an independent set of size $2d-2$;
whereas if $\gS$ is the complete $d$-regular bipartite graph then
$\ell(\gS)=d$.)

\medskip

Our main result is the following. Recall that $H(x)=-x\log x
-(1-x)\log (1-x)$ is the usual binary entropy function.
\begin{thm} \label{thm-main}
Fix $\rho>0$ satisfying $H(\rho)+\rho < 1$. There are constants
$d_0, C_1, C_1', C_2>0$ all depending on $\rho$ such that if $\gS$
is a $d$-regular bipartite graph on $N$ vertices with bipartite
expansion $\gd$ and locality $\ell>0$ satisfying
\begin{equation} \label{inq-bounds.on.gd}
\gd \geq \max\left\{\frac{C_1\log^3 d}{d},\frac{C_1'\log
d}{\ell}\right\}
\end{equation}
and with $d \geq d_0$ and if $\cM(\gS)$ is an ergodic $\rho$-local
Markov chain on state space $\cC_3(\gS)$ with stationary
distribution $\pi_3(\gS)$ then
$$
\tau_{\cM(\gS)} \geq \exp_2\left\{\frac{C_2N\gd}{\log d}\right\}.
$$
\end{thm}
Note that for all $\rho \leq .22$ we have $H(\rho)+\rho < 1$. Here
and throughout we use ``$\log$'' for $\log_2$ and write $\exp_2 x$
for $2^x$.

\begin{remark}
The second inequality in (\ref{inq-bounds.on.gd}) implies $\ell \geq
\Omega(\log d/\gd)$. This condition appears in the derivation of
(\ref{done.for.small}), where it is only used in the weaker form
$\ell=\omega(1)$ (which follows since $\gd \leq 1$). It is used in a
more essential way in the derivation of (\ref{using.iso.2}) where it
serves to limit, somewhat artificially, the number of $3$-colourings
of a bipartite graph with a given pre-image of $0$. We expect that
Theorem \ref{thm-main} should remain true with the second inequality
in (\ref{inq-bounds.on.gd}) removed.
\end{remark}

We now return to Glauber dynamics. This changes the colour of at
most one vertex at each step, and so (as long as the underlying
graph has at least five vertices) it is a $\rho$-local chain for
$\rho=.2$. Fixing $\rho$ to this value, all of the constants in
Theorem \ref{thm-main} become absolute, and we have the following
corollary.
\begin{cor} \label{cor-1}
There are constants $d_0, C_1, C_1', C_2>0$ such that if $\gS$
satisfies the conditions of Theorem \ref{thm-main} and if the
Glauber dynamics chain $\cM_3(\gS)$ is ergodic then
$$
\tau_{\cM_3(\gS)} \geq \exp_2\left\{\frac{C_2N\gd}{\log
d}\right\}.
$$\end{cor}

Let us apply Theorem \ref{thm-main} to the case $\gS=Q_d$, the
$d$-dimensional Hamming cube. This is the graph on vertex set
$\{0,1\}^d$ with $x\sim y$ iff $x$ and $y$ differ on exactly one
coordinate. For $d\geq 2$ we have $\ell(Q_d)=d$ (the graph induced
by the union of the neighbourhoods of adjacent vertices is a perfect
matching, so all independent sets are of size at most $d$), and
$\gd(Q_d) = \Omega(1/\sqrt{d})$ (see {\em e.g.} \cite[Lemma
1.3]{KorshunovSapozhenko}). So the following is an immediate
corollary of Theorem \ref{thm-main}.
\begin{cor} \label{cor-cube.1}
Fix $\rho>0$ satisfying $H(\rho)+\rho < 1$. There is a constant
$C=C(\rho)>0$ such that for all $d \geq 2$, if $\cM(Q_d)$ is an
ergodic $\rho$-local Markov chain on state space $\cC_3(Q_d)$ with
stationary distribution $\pi_3(Q_d)$ then
$$
\tau_{\cM} \geq \exp_2\left\{\frac{C2^d}{\sqrt{d}\log d}\right\}.
$$
\end{cor}
In particular this result applies to the Glauber dynamics chain,
although in this case it is not necessary to hypothesize ergodicity.
\begin{cor} \label{cor-cube.2}
There is a constant $C>0$ such that for all $d \geq 2$,
$$
\tau_{\cM_3(Q_d)} \geq \exp_2\left\{\frac{C2^d}{\sqrt{d}\log
d}\right\}.
$$
\end{cor}

\noindent {\em Proof: }In the presence of Corollary
\ref{cor-cube.1}, it suffices to show that the chain $\cM_3(Q_d)$ is
ergodic. We will show that if $\chi_1$ is a $3$-colouring of $Q_d$
with $\chi_1(v_0)=0$ for some $v_0 \in \cE$ then there is a sequence
of steps in the Glauber dynamics chain that takes $\chi_1$ to a
$2$-colouring $\chi_2$ of $Q_d$ with $\chi_2(v)=0$ for all $v \in
\cE$. This suffices, since it is clear that any one of the six
$2$-colourings of $Q_d$ can be reached from any other via steps in
the chain.

We make use of a correspondence between proper $3$-colourings of
$Q_d$ and homomorphisms from $Q_d$ to $\Z$ that send $v_0$ to $0$.
Formally, set
$$
\cF^{v_0}=\{f:V \rightarrow \Z:f(v_0)=0~\mbox{and}~x\sim y
\Rightarrow |f(x)-f(y)|=1\}.
$$
(This set was introduced in \cite{BenjaminiHaggstromMossel} and
further studied in \cite{Galvin,Kahn}.) Then, as observed by Randall
\cite{Randall2}, there is a bijection from $\cF^{v_0}$ to
$\cC_3^{v_0}:=\{\chi \in \cC_3:\chi(v_0)=0\}$ given by $f
\longrightarrow \Phi(f)$ where $\Phi(f)(v)=i$ iff $f(v) \equiv i$
(mod $3$). Before verifying that this is indeed a bijection, we use
the correspondence to establish the corollary.

For $f \in \cF^{v_0}$ set $R(f)=\{f(v):v \in V\}$. Now consider
$\chi_1 \in \cC_3^{v_0}$. If $|R(\Phi^{-1}(\chi_1))|=2$, then we
may take $\chi_2=\chi_1$ and we are done. If
$|R(\Phi^{-1}(\chi_1))|=k>2$, then it suffices to exhibit a
sequence of steps in the chain that takes $\chi_1$ to some $\chi_3
\in \cC_3^{v_0}$ with $|R(\Phi^{-1}(\chi_3))|=k-1$.

Without loss of generality we may assume that $\Phi^{-1}(\chi_1)$
takes on some strictly positive values. Let $\ell$ be the largest
such value, and let $v \in V$ be any vertex satisfying
$\Phi^{-1}(\chi_1)(v)=\ell$. Note that $\ell-2 \in
R(\Phi^{-1}(\chi_1))$. Let $f:V\rightarrow \Z$ be the function
that agrees with $\Phi^{-1}(\chi_1)$ off $v$ and satisfies
$f(v)=\ell-2$. Since $\Phi^{-1}(\chi_1)(y)=\ell-1$ for all $y \in
N(v)$ and $v \neq v_0$ we have that $f \in \cF^{v_0}$ and $\Phi(f)
\in \cC_3^{v_0}$ and that the Glauber dynamics chain permits a
move from $\chi_1$ to $\Phi(f)$. But we also have that $R(f)
\subseteq R(\Phi^{-1}(\chi_1))$ and $|\{v \in V:f(v)=\ell\}|<|\{v
\in V:\Phi^{-1}(\chi_1)(v)=\ell\}|$, so that by repeating the
above described procedure $m$ more times (where $m=|\{v \in
V:f(v)=\ell\}|$) we arrive at the desired $\chi_3$.

It remains to verify that $\Phi$ is a bijection. That it is
injective is clear. To see that it is surjective, consider $\chi'
\in \cC_3^{v_0}$. We shall construct from $\chi'$ an $f \in
\cF^{v_0}$ with $\Phi(f)=\chi'$ by setting $f(v_0)=0$ and then
extending $f$ level by level, where the $k^{th}$ level of $Q_d$
($k=0,\ldots,d)$ is ${\cal L}_k:=\{v \in V:dist(v,v_0)=k\}$ (here
we are using $dist(\cdot,\cdot)$ for the usual graph distance).
Note that for $v \in {\cal L}_k$, $N(v) \subseteq {\cal L}_{k-1}
\cup {\cal L}_{k+1}$ and that for $f \in \cF^{v_0}$ the values
that $f$ takes on ${\cal L}_k$ must all have the same parity.

So suppose we have specified $f$ up to ${\cal L}_k$ for some $0
\leq k \leq d-1$. Consider $v \in {\cal L}_{k+1}$. If $f$ is
constant on $N(v) \cap {\cal L}_k$ then (since the construction of
$f$ has succeeded up to ${\cal L}_k$) we also have that $\chi'$ is
constant on $N(v) \cap {\cal L}_k$ with $\chi'(y) \equiv f(y)$
(mod $3$) for all $y \in N(v) \cap {\cal L}_k$. In this case we
choose $f(v)$ such that $|f(v)-f(y)|=1$ for all $y \in N(v) \cap
{\cal L}_k$ and $\chi'(v) \equiv f(v)$ (mod $3$).

If $f$ is not constant on $N(v) \cap {\cal L}_k$, then we claim
that there is some $\ell \in \Z$ such that $f$ takes on only the
values $\ell$ and $\ell + 2$ on $N(v) \cap {\cal L}_k$. For if
not, then we have $y_1,y_2 \in N(v) \cap {\cal L}_k$ with
$|f(y_1)-f(y_2)|\geq 4$. But by the structure of $Q_d$ there must
be $v' \in {\cal L}_{k-1}$ with $v' \sim y_1$ and $v' \sim y_2$,
which forces $|f(y_1)-f(y_2)|\leq 2$. This contradiction
establishes the two-value claim. We now set $f(v)=\ell+1$,
allowing the construction to continue. Since ${\cal L}_{k+1}$ is
an independent set in $Q_d$, we may repeat the above-described
procedure on each vertex of ${\cal L}_{k+1}$ independently, thus
extending the construction of $f$ to all of ${\cal L}_{k+1}$. \qed

\begin{remark}
Glauber dynamics for $q$-colourings of $Q_d$ is not in general
ergodic for $3 < q < \gD(Q_d)+1$. Indeed, it is straightforward to
construct a $4$-colouring $\chi$ of $Q_3$ which is frozen in the
sense that $P_4(\chi,\chi')=0$ for all $\chi' \neq \chi$; one simply
assigns the colours $0$, $1$, $2$ and $3$ to a particular vertex and
its three neighbours and then extend to a colouring of the whole of
$Q_3$ according to the rule that on each face ($4$-cycle) of $Q_3$
all of the colours $0$, $1$, $2$ and $3$ must appear.
\end{remark}

\begin{remark}
While this paper was under review, Galvin and Randall
\cite{GalvinRandall} used methods different to those of the present
work to extend Corollary \ref{cor-cube.2} to the discrete torus
$T_{L,d}$, the graph on vertex set $\{0, \ldots, L-1\}^d$ in which
two strings are adjacent if they differ on exactly one coordinate
and differ by $1~(\rm{mod}~L)$ on that coordinate. The main result
of \cite{GalvinRandall} is that for $L \geq 4$ even and $d$ large,
the Glauber dynamics chain $\cM_3$ on $\cC_3(T_{L,d})$ satisfies
$\tau_{\cM_3} \geq \exp \{L^{d-1}/(d^4\log^2 L)\}$.
\end{remark}

We prove Theorem \ref{thm-main} via a well-known conductance
argument (introduced in \cite{JerrumSinclair}). A particularly
useful form of the argument was given by Dyer, Frieze and Jerrum
\cite{DyerFriezeJerrum}. Let $\cM$ be an ergodic Markov chain on
state space $\gO$ with transition probabilities $P$ and stationary
distribution $\pi$. Let $A \subseteq \gO$ and $M \subseteq \gO
\setminus A$ satisfy $\pi(A) \leq 1/2$ and $\go_1 \in A, \go_2 \in
\gO \setminus (A \cup M) \Rightarrow P(\go_1, \go_2) =0$. Then from
\cite{DyerFriezeJerrum} we have
$$
\tau_\cM \geq \frac{\pi(A)}{8\pi(M)}.
$$
We may think of $M$ as a bottleneck set through which any run of the
chain must pass in order to mix; if the bottleneck has small
measure, then the mixing time is high.

Now let us return to the setup of Theorem \ref{thm-main}. Set
$$
\cC_3^{b,\rho,0} = \cC_3^{b,\rho,0}(\gS) = \{\chi \in
\cC_3:\left||\chi^{-1}(0)\cap \cE|-|\chi^{-1}(0)\cap \cO|\right|
\leq \rho N/2\}
$$
($\cC_3^{b,\rho,0}$ is the set of $3$-colourings that are balanced
with respect to $0$) and
$$
\cC_3^{\cE,\rho,0} = \cC_3^{\cE,\rho}(\gS) = \{\chi \in
\cC_3:|\chi^{-1}(0)\cap \cE|>|\chi^{-1}(0)\cap \cO| + \rho N/2\}.
$$
We may assume without loss of generality that
$\pi_3(\cC_3^{\cE,\rho,0}) \leq 1/2$. Notice that since $\cM$
changes the colour of at most $\rho N$ vertices in each step, we
have that if $\chi_1 \in \cC_3^{\cE,\rho,0}$ and $\chi_2 \in \cC_3
\setminus (\cC_3^{\cE,\rho,0} \cup \cC_3^{b,\rho,0})$ then
$P_\cM(\chi_1,\chi_2) = 0$. We therefore have
$$
\tau_\cM \geq
\frac{\pi_3(\cC_3^{\cE,\rho,0})}{8\pi_3(\cC_3^{b,\rho})} \geq
\frac{2^{N/2}}{8|\cC_3^{b,\rho,0}|},
$$
the second inequality coming from the trivial lower bound
$|\cC_3^{\cE,\rho,0}|\geq 2^{N/2}$ (consider those $\chi$ with
$\chi(v)=0$ for all $v \in \cE$). Theorem \ref{thm-main} thus
follows from the following theorem, whose proof will be the main
business of this paper.
\begin{thm} \label{thm-main2}
Fix $\rho>0$ satisfying $H(\rho)+\rho < 1$. There are constants
$d_0, C_1, C_1', C_2>0$ all depending on $\rho$ such that if $\gS$
is a $d$-regular bipartite graph on $N$ vertices with bipartite
expansion $\gd$ and locality $\ell$ satisfying
(\ref{inq-bounds.on.gd}) and with $d \geq d_0$ then
$$
|\cC_3^{b,\rho,0}| \leq
\exp_2\left\{\frac{N}{2}\left(1-\frac{C_2\gd}{\log
d}\right)\right\}.
$$
\end{thm}

Theorem \ref{thm-main2} says more about the structure of $\cC_3$
than just that the dynamics mixes slowly. From it, we can infer that
for $\gS$ and $\rho$ satisfying the conditions of the theorem,
$\cC_3$ breaks naturally into six sets in such a way that once a
$\rho$-local chain enters one of these dominant sets, it tends to
remain there for an exponential time. These sets are characterized
by a predominance of one (of three) colours on one (of two)
partition classes. Indeed, defining $\cC_3^{b,\rho,1}$ and
$\cC_3^{b,\rho,2}$ by analogy with $\cC_3^{b,\rho,0}$ and setting
$\cR_3 = \cC_3 \setminus \cup_{i=0}^2 \cC_3^{b,\rho,i}$, we may
partition $\cR_3$ into six pieces by
$$
\cR_3 = \cup_{(x,y,z) \in
\{\cE,\cO\}^3\setminus\{(\cE,\cE,\cE),(\cO,\cO,\cO)\}}
\cR_3^{(x,y,z)}
$$
where $\cR_3^{(x,y,z)} = \{\chi \in \cR_3 : \chi \in
\cC_3^{x,\rho,0} \cap \cC_3^{y,\rho,1} \cap \cC_3^{z,\rho,2}\}$. If
a $\rho$-local chain leaves $\cR_3^{(x,y,z)}$ (for any $(x,y,z)$) it
must enter $\cup_{i=0}^2 \cC_3^{b,\rho,i}$ which, by Theorem
\ref{thm-main} and a union bound, has exponentially small measure.

Before turning to the proof of Theorem \ref{thm-main2} we pause to
give a pleasing combinatorial corollary in the special case
$\gS=Q_d$.
\begin{cor} \label{cor-imbalance in cube}
Fix $\rho$ satisfying $H(\rho)+\rho<1$. There is a constant
$C=C(\rho)>0$ such that for all $d \geq 2$, if $\chi$ is a uniformly
chosen $3$-colouring of $Q_d$ then
$$
{\mathbb P}\left(\exists i:\left|\frac{|\chi^{-1}(i)\cap
\cE|}{|\cE|}-\frac{|\chi^{-1}(i)\cap \cO|}{|\cO|}\right| \leq
\rho\right) \leq \exp_2\left\{-\frac{C2^d}{\sqrt{d}\log d}\right\}.
$$
\end{cor}
In other words, the typical $3$-colouring of $Q_d$ exhibits strong
$\cE/\cO$ imbalance on all colours.

\medskip

\noindent {\em Proof of Corollary \ref{cor-imbalance in cube}: }As
previously observed, $\ell(Q_d)=d$ and $\gd(Q_d)\leq
\Omega(1/\sqrt{d})$, so (\ref{inq-bounds.on.gd}) is satisfied for
large enough $d$. It follows that there is a $C'(\rho)$ such that
for large enough $d$ and for each $i=0,1,2$,
$$
\left|\{\chi \in \cC_3(Q_d): \left||\chi^{-1}(i)\cap
\cE|-|\chi^{-1}(i)\cap \cO|\right| \leq \rho 2^{d-1}\}\right| \leq
\exp_2\left\{2^{d-1}-\frac{C'2^d}{\sqrt{d}\log d}\right\}.
$$
Using $2^{2^{d-1}}$ as a lower bound on $|\cC_3(Q_d)|$ (consider
those colourings for which $\chi^{-1}(0)=\cE$) we obtain
$$
{\mathbb P}\left(\left||\chi^{-1}(i)\cap \cE|-|\chi^{-1}(i)\cap
\cO|\right| \leq \rho 2^{d-1}\right) \leq
\exp_2\left\{-\frac{C'2^d}{\sqrt{d}\log d}\right\}
$$
for $\chi$ chosen uniformly from $\cC_3(Q_d)$. The stated bound
follows for large $d$ (with a constant $C''$ slightly larger than
$C'$) via a union bound and the fact that $|\cE|=|\cO|=2^{d-1}$; we
may obtain the bound for all $d$ by appropriately modifying the
constant $C''$. \qed

\section{Proof of Theorem \ref{thm-main2}} \label{sec-proof.of.mn.thm}

\subsection{Overview of the proof}

In this section we give an informal overview of the proof of Theorem
\ref{thm-main2}.

We bound the number of balanced $3$-colourings by bounding, for
each pair $E \subseteq \cE$, $O \subseteq \cO$ with $||E|-|O||
\leq \rho N/2$ and $E \not \sim O$ (that is, with no edge in $\gS$
joining $E$ and $O$), the number of $3$-colourings of $\gS$ in
which $E \cup O$ is the pre-image of $0$. We then sum over all
choices of $E$ and $O$.

How many ways are there to $3$-colour $\gS$ given that $E \cup O$ is
the pre-image of $0$? Write $I(E)$ for the set of vertices in $N(E)$
all of whose neighbours are in $E$, and $I(O)$ for the neighbours of
$O$ all of whose neighbours are in $O$. There are two choices for
each vertex in $I(E)$ and two for each vertex in $I(O)$, as well as
two choices for each component in the graph obtained from $\gS$ by
removing $E$, $O$, $I(E)$ and $I(O)$ (each such component is a
connected bipartite graph), all choices independent. The first step
in the proof is an easy graph theory lemma that shows that the
contribution from components in $\gS - (E \cup O \cup I(E) \cup
I(O))$ is negligible. (This step uses the locality of $\gS$ in an
essential way.) This reduces the problem of bounding the number of
balanced $3$-colorings to the problem of estimating a sum of the
form
\begin{equation}
\sum_{E, O: E \not \sim O} 2^{|I(E)|+|I(O)|}. \label{sum-inf}
\end{equation}
When $|E|$ and $|O|$ are both small (less than $cN$ for a suitably
small constant $c$) a naive count suffices to give an appropriate
bound. For larger $E$ and $O$, we must work harder. We partition
the set of pairs $(E,O)$ according to the parameters $a=|[E]|$,
$g=|N(E)|$, $b=|I(E)|$, $h=|N(I(E))|$, $b'=|I(O)|$ and
$h'=|N(I(O))|$. Within each class, each pair gives the same
contribution ($2^{b+b'}$) to the sum in (\ref{sum-inf}). The main
point of the proof is an estimate on the size of
$\cH=\{(E,O):(E,O)$ has parameters $a$, $g$, $b$, $h$, $b'$ and
$h'\}$ of the form
\begin{equation} \label{main.est-inf}
|\cH| \leq \exp_2\left\{\frac{N}{2}-b-b'-\frac{cN\gd}{\log
d}\right\}
\end{equation}
for sufficiently large $d=d(\rho)$ and suitable $c=c(\rho)$. The
proof is completed by invoking (\ref{main.est-inf}) and summing
over all choices of $a$, $g$, {\em et cetera}.

The proof of (\ref{main.est-inf}) involves the idea of
approximation. To bound $|\cH|$, we produce a small set $\cU$ with
the properties that each $(E,O) \in \cH$ is approximated (in an
appropriate sense) by some $U \in \cU$, and for each $U \in \cU$,
the number of $(E,O) \in \cH$ that could possibly be approximated by
$U$ is small. (Each $U \in {\cal U}$ will consist of six parts; one
each approximating $E$, $N(E)$, $I(E)$, $N(I(E))$, $I(O)$ and
$N(I(O))$.) The product of the bound on $|\cU|$ and the bound on the
number of those $(E,O) \in \cH$ that may be approximated by any $U$
is then a bound on $|\cH|$.

The main inspiration for our approximation scheme is the work of
A. Sapozhenko, who, in \cite{Sapozhenko2}, gave a relatively
simple derivation for the asymptotics of the number of independent
sets in $Q_d$, earlier derived in a more involved way in
\cite{KorshunovSapozhenko}. We produce the set $\cU$ by appealing
to a lemma from \cite{GalvinTetali} where a similar approximation
scheme was used to show that the mixing time of Glauber dynamics
for the hard-core model on $Q_d$ with activity $\gl$ is
(essentially) exponential in $2^d$ for large enough $\lambda$. The
proof that each $U \in \cU$ approximates only a small number of
$(E,O) \in \cH$ is a modification of a similar proof from
\cite{Galvin} in which it is shown that a uniformly chosen
homomorphism from $Q_d$ to $\Z$ almost surely takes on at most $5$
values, and also that the number of proper $3$-colourings of $Q_d$
is asymptotic to $2e2^{2^{d-1}}$ as $d$ goes to infinity.

\subsection{The proof}

We begin by establishing some more notation. From now on, we write
$M$ for $N/2$. For $A \subseteq \cE$ and $B \subseteq \cO$ write
$A \not \sim B$ if for all $x \in A$ and $y \in B$, $x \not \sim
y$ (this is equivalent to both $N(A) \cap B = \emptyset$ and $N(B)
\cap A = \emptyset$). For $S \subseteq V$ write $d_S(u)$ for
$|N(u) \cap S|$ and $\comp(S)$ for the number of components of the
subgraph induced by $S$. Finally for $T \subseteq \cE$ (or $\cO$)
set
$$
I(T) = \{x \in N(T) : N(x) \subseteq T\}~(=\{x \in V : N(x)
\subseteq T\}).
$$
We think of $I(T)$ as an internal closure of $N(T)$. Note that for
all $T \subseteq \cE$ (or $\cO$), $[I(T)]=I(T)$, $I(T) \subseteq
N(T)$ and $N(I(T)) \subseteq T$.

For $\chi \in \cC_3^{b,\rho,0}$ set
$$
E=\chi^{-1}(0)\cap \cE,
$$
$$
O=\chi^{-1}(0)\cap \cO,
$$
$$
I=I(E),
$$
$$
J=I(O)
$$
and
$$
R=V\setminus (E\cup O\cup I \cup J).
$$
We assume the convention that whenever $E$ and $O$ have been
specified, $I$, $J$ and $R$ will be used as shorthand for $I(E)$,
$I(O)$ and $V\setminus (E\cup O\cup I \cup J)$.

For $E \subseteq \cE$ and $O \subseteq \cO$ set
$$
\cC_3(E,O)=\{\chi \in \cC_3 : E(\chi)=E,~O(\chi)=O\}.
$$
Note that $\cC_3(E,O) \neq \emptyset$ iff $E \not \sim O$. For
$\cC_3(E,O) \neq \emptyset$ we have
$$
|\cC_3(E,O)|=2^{|I|+|J|+\comp(R)}.
$$
To see this, note that once we have specified that the set of
vertices coloured $0$ is $E \cup O$, we have a free choice between
$1$ and $2$ for the colour at $x \in I \cup J$, with each choice
independent. This accounts for the factor $2^{|I|+|J|}$. The
subgraph induced by $R$ breaks into $\comp(R)$ components, each of
which is bipartite and may be coloured in exactly two ways using
the colours $1$ and $2$. This accounts for the factor
$2^{\comp(R)}$. We therefore have
$$
|\cC_3^{b,\rho,0}| = \sum_{E \subseteq \cE,~O \subseteq \cO:
\atop{||E|-|O|| \leq \rho M,~E\not \sim O}} 2^{|I|+|J|+\comp(R)}.
$$

A key observation is the following.

\begin{prop} \label{inq-number.of.remaining.comps}
For $E \not \sim O$, $\comp(R) \leq 2M/\ell$.
\end{prop}

\noindent {\em Proof: }Let $C$ be a component of $V \setminus (E
\cup O)$. If $C=\{v\}$ consists of a single vertex, then
(depending on the parity of $v$) we have either $N(v) \subseteq E$
or $N(v) \subseteq O$ and so $v \in I(E) \cup I(O)$. Otherwise,
let $vw$ be an edge of $C$. We have
$$
|(N(v) \cup N(w)) \cap (E
\cup O)| \leq 2d-\ell
$$
(recall that $E \cup O = \chi^{-1}(0)$ is
an independent set), and so $|C|\geq \ell$. The result follows.
\qed

\medskip

We now decompose $\cC_3^{b,\rho,0}$ into four pieces. Set
\begin{equation} \label{def-alpha}
\ga = \sup \left\{\ga ' \in \left[0,\frac{1}{2}-\rho\right] : 2\ga'+
\rho + H(\ga')+H(\rho+ \ga') \leq \frac{1}{2}\left(1 + \rho +
H(\rho)\right)\right\}.
\end{equation}
Since $H(\rho)+\rho<1$, $\ga$ is a strictly positive constant
depending on $\rho$. Set
$$
\cC_3^{b,\rho,0}(triv,\cE) = \{\chi \in \cC_3^{b,\rho,0} : |E| \leq
\ga M, |E| \leq |O|\}
$$
and define $\cC_3^{b,\rho,0}(triv,\cO)$ analogously. Set
$$
\cC_3^{b,\rho,0}(nt,\cE) = \{\chi \in \cC_3^{b,\rho,0} \setminus
(\cC_3^{b,\rho,0}(triv,\cE) \cup \cC_3^{b,\rho,0}(triv,\cO)) : E
~\mbox{small}\}
$$
(recall that $E$ is small if $|[E]| \leq M/2$) and define
$\cC_3^{b,\rho,0}(nt,\cO)$ similarly. Since $\gS$ has a perfect
matching, it is easy to see that for $\chi \in \cC_3$ at least one
of $|E| \leq M/2$, $|O| \leq M/2$ holds; moreover, it is
straightforward to check that at least one of $|[E]|\leq M/2$,
$|[O]|\leq M/2$ holds also; that is, that at least one of $E$, $O$
is small, and so
$$
\cC_3^{b,\rho,0} = \cC_3^{b,\rho}(triv,\cE) \cup
\cC_3^{b,\rho,0}(triv,\cO) \cup \cC_3^{b,\rho}(nt,\cE) \cup
\cC_3^{b,\rho,0}(nt,\cO).
$$

In what follows we make extensive use of a result concerning the
sums of binomial coefficients which follows from the Chernoff
bounds \cite{Chernoff} (see also \cite{Bollobas3}, p.11):
\begin{equation} \label{inq-binomial}
\sum_{i=0}^{[\gb M]}{M\choose i} \leq 2^{H(\gb)M}~~~~~\mbox{for $\gb
\leq \frac{1}{2}$}.
\end{equation}
Also, since $H(x)\leq 2x\log 1/x$ for $x \leq e^{-1}$,
\begin{equation} \label{inq-binomial.2}
\sum_{i=0}^{[\gb M]}{M\choose i} \leq 2^{2\gb
M\log(1/\gb)}~~~~~\mbox{for $\gb \leq e^{-1}$}.
\end{equation}

We begin by bounding $|\cC_3^{b,\rho,0}(triv, \cE)|$. Noting that
$|I|\leq |E|$ and $|J| \leq |O|$ always (this follows from
(\ref{obs-bipartite,regular})) we have
\begin{eqnarray}
|\cC_3^{b,\rho,0}(triv,\cE)| & = & \sum_{E \subseteq \cE,~O
\subseteq \cO: \atop{|E| \leq \ga M,~|O|\leq (\ga + \rho) M}}
2^{|E|+|O|+\comp(R)} \nonumber \\
& \leq & \exp_2\left\{\frac{2M}{\ell}\right\} \sum_{E \subseteq
\cE,~O \subseteq \cO: \atop{|E| \leq \ga M,~|O|\leq (\ga + \rho) M}}
2^{|E|+|O|} \label{using_prop_2.1} \\
& \leq & \exp_2\left\{M\left(2\ga + \rho +
\frac{2}{\ell}\right)\right\} \sum_{i=0}^{[\ga M]} {M
\choose i} \sum_{i=0}^{[(\ga + \rho)M]} {M \choose i} \nonumber \\
& \leq & \exp_2\left\{M\left(2\ga + \rho + H(\ga) + H(\rho + \ga)
+ \frac{2}{\ell}\right)\right\}
\label{using.entropy} \\
& \leq & \exp_2\left\{\frac{N}{2}\left(1-\frac{\gd}{\log
d}\right)\right\} \label{done.for.small}
\end{eqnarray}
for sufficiently large $d=d(\rho)$. In (\ref{using_prop_2.1}) we
have used Proposition \ref{inq-number.of.remaining.comps}. In
(\ref{using.entropy}) we use (\ref{inq-binomial}) while in
(\ref{done.for.small}) we use (\ref{def-alpha}) to obtain
$$
2\ga + \rho + H(\ga) + H(\rho + \ga) + \frac{2}{\ell} \leq
1-\varepsilon + \frac{2}{\ell}
$$
for some $\varepsilon=\varepsilon(\ga)$, and then use the second
inequality in (\ref{inq-bounds.on.gd}) (in the weak form that
$\ell=\omega(1)$) to obtain
$$
1-\varepsilon + \frac{2}{\ell} \leq 1-\frac{\gd}{\log d}
$$
(note that $\gd/\log d =o(1)$). Similarly, we have
\begin{equation} \label{done.for.small.2}
|\cC_3^{b,\rho,0}(triv,\cO)| \leq
\exp_2\left\{\frac{N}{2}\left(1-\frac{\gd}{\log d}\right)\right\}
\end{equation}
for suitable $d$.

Next we turn to $\cC_3^{b,\rho,0}(nt,\cE)$ and
$\cC_3^{b,\rho,0}(nt,\cO)$. Without loss of generality we may assume
$|\cC_3^{b,\rho,0}(nt,\cE)| \leq |\cC_3^{b,\rho,0}(nt,\cO)|$.
Bearing (\ref{done.for.small}) and (\ref{done.for.small.2}) in mind,
Theorem \ref{thm-main2} now follows from
\begin{equation} \label{to.prove}
|\cC_3^{b,\rho,0}(nt,\cE)| \leq
\exp_2\left\{\frac{N}{2}\left(1-\frac{c\gd}{\log d}\right)\right\}
\end{equation}
for some constant $c=c(\rho)$.

For integers $a, g, b, h, b'$ and $h'$, set
$$
\cH(a,g,b,h,b',h') = \left\{(E,O): {E \subseteq \cE,~O \subseteq
\cO,~E \not \sim O,~|[E]|=a,~|N(E)|=g,
\atop{|I|=b,~|N(I)|=h,~|J|=b',~|N(J)|=h'}}\right\}.
$$
Our main lemma is the following ({\em cf.} \cite[Theorem
2.1]{GalvinTetali}).

\begin{lemma} \label{lem-main}
For each $\gb_1, \gb_2 > 0$, there are constants $d_0, c>0$
depending on both $\gb_1$ and $\gb_2$ such that the following holds.
If $G$ is a $d$-regular bipartite graph with bipartite expansion
$\gd \geq d^{-\gb_1}$ and $d\geq d_0$ and if $a$ satisfies $\gb_2 M
\leq a \leq M/2$, then for any $g,b,h,b'$ and $h'$ we have
$$
|\cH(a,g,b,h,b',h')| \leq \exp_2\left\{M\left(1+ \frac{15\log^2
d}{d}\right)-b-b'-\frac{c\gd g}{\log d}\right\}.
$$
\end{lemma}

For $G=\gS$ we may take $\gb_1=2$ (say). Note that for each $(E,O)
\in \cC_3^{b,\rho,0}(nt,\cE)$ with $|I|=b$ and $|J|=b'$, $(E,O) \in
\cH(a,g,b,h,b',h')$ for some $a$, $g$, $h$ and $h'$ with $\ga M \leq
a \leq M/2$. With the steps justified below, we therefore have
\begin{eqnarray}
|\cC_3^{b,\rho,0}(nt,\cE)| & \leq &
     \exp_2\left\{\frac{2M}{\ell}\right\}\sum_{a,g,b,h,b',h': \atop{\ga M \leq a \leq M/2}}
     |\cH(a,g,b,h,b',h')| 2^{b+b'} \nonumber \\
& \leq &
     \exp_2\left\{M\left(1+\frac{2}{\ell} + \frac{15\log^2 d}{d}\right)\right\}
     \sum_{a,g,b,h,b',h': \atop{\ga M \leq a \leq M/2}}
     \exp_2\left\{-\frac{c\gd g}{\log d}\right\} \label{using.main.lem} \\
& \leq &
     \exp_2\{M\left(1+\frac{2}{\ell} + \frac{21\log^2 d}{d}\right)
     \max_{\ga M \leq a \leq M/2 \atop{a\leq g}}
     \left\{\exp_2\left\{-\frac{c\gd g}{\log d}\right\}\right\} \label{d.large} \\
& \leq &
     \exp_2\left\{M\left(1+\frac{2}{\ell} + \frac{21\log^2 d}{d}
     - \frac{c\ga\gd}{\log d}\right)\right\} \label{using.iso.1} \\
& \leq &
     \exp_2\left\{\frac{N}{2}\left(1-\frac{c'\gd}{\log
     d}\right)\right\} \label{using.iso.2}
\end{eqnarray}
verifying (\ref{to.prove}) and completing the proof of Theorem
\ref{thm-main2}. The main point, (\ref{using.main.lem}), is an
application of Lemma \ref{lem-main}. Here the constant $c$ depends
on $\ga$ and therefore on $\rho$. In (\ref{d.large}) we use that $M
\leq \exp\{M \log^2 d/d\}$ for all $d$. In (\ref{using.iso.1}) we
have chosen $g = \ga M$ to maximize the exponent. Finally in
(\ref{using.iso.2}) we may (for example) take $C_1 = 43/(\ga c)$ and
$C_1' = 4/(\ga c)$, and use both inequalities in
(\ref{inq-bounds.on.gd}). The final constant $c'$ depends only on
$c$ and $\ga$ and therefore only on $\rho$, as claimed.

\medskip

To prove Lemma \ref{lem-main}, we use a notion of approximation
introduced in \cite{Sapozhenko}. An {\em approximation} for $A
\subseteq \cE$ is a pair $(F,S) \subseteq \cO \times \cE$ satisfying
$$
F \subseteq N(A),~~S \supseteq [A],
$$
$$
d_F(u) \geq d-\sqrt{d}~~~\forall u \in S
$$
and
$$
d_{\cE \setminus S}(v) \geq d-\sqrt{d}~~~\forall v \in \cO \setminus
F.
$$
For $A \subseteq \cO$ we make the analogous definition.

The following lemma is from \cite{GalvinTetali} (a combination of
Lemmata 3.2 and 3.3). We use the shorthand ${t \choose \leq k}$ for
$\sum_{0 \leq i \leq k} {t \choose i}$.

\begin{lemma} \label{lem-basic.approx}
Let $G$ be a $d$-regular bipartite graph with $2M$ vertices. For
each $a$ and $g$ set
$$
\cA(a,g)=\{A \subseteq \cE: |[A]|=a, |N(A)|=g\}.
$$
There is a family $\cW=\cW(a, g) \subseteq 2^\cO \times 2^\cE$ with
$$
|\cW| \leq {M \choose \leq \frac{2g\log d}{d}}{2g\log d \choose \leq
\frac{2g}{d}}{2d^3g\log d \choose \leq
\frac{2(g-a)}{\sqrt{d}}}{2g\log d \choose \leq
(g-a)\frac{\sqrt{d}}{(d-\sqrt{d})}}
$$
such that every $A \in \cA(a,g)$ has an approximation in $\cW$. The
analogous result holds with $\cO$ replacing $\cE$ in the definition
of $\cA(a,g)$.
\end{lemma}

\begin{remark} \label{rem-specific.approx}
If $a$ and $g$ satisfy $g-a \geq d^{-\gb}g$ for some constant $\gb$
then using (\ref{inq-binomial.2}) the bound on $|\cW|$ from Lemma
\ref{lem-basic.approx} may be rewritten as
$$
|\cW| \leq \exp_2 \left\{\frac{5M\log^2d}{d} +
\frac{(6\gb+17)(g-a)\log d}{\sqrt{d}}\right\}
$$
as long as $d$ is sufficiently large (as a function of $\gb$).
\end{remark}

\medskip

Say that a sextuple $(F,S,P,Q,P',Q') \subseteq \cO \times \cE \times
\cE \times \cO \times \cO \times \cE$ is an {\em approximation} for
$(E,O) \in \cH(a,g,b,h,b',h')$ if $(F,S)$ is an approximation for
$E$, $(P,Q)$ is an approximation for $I$ and $(P',Q')$ is an
approximation for $J$.

\begin{lemma} \label{lem-adv.approx}
Let $G$, $a$, $g$, $b$, $h$, $b'$ and $h'$ be as in Lemma
\ref{lem-main}. There are constants $c_1=c_1(\gb_1)>0$ and
$c_2=c_2(\gb_1,\gb_2)>0$ and a family $\cX=\cX(a,g,b,h,b',h')
\subseteq 2^\cO \times 2^\cE \times 2^\cE \times 2^\cO \times
2^\cO \times 2^\cE$ with
$$
|\cX| \leq \exp_2\left\{\frac{15M\log^2 d}{d}+ \frac{c_1(g-a)\log
d}{\sqrt{d}} + \frac{c_1(h-b)\log d}{\sqrt{d}} +
\frac{c_2(h'-b')\log d}{\sqrt{d}} \right\}
$$
such that every $(E,O) \in \cH(a,g,b,h,b',h')$ has an approximation
in $\cX$.
\end{lemma}

\noindent {\em Proof: }We apply Lemma \ref{lem-basic.approx} (in
the form given in Remark \ref{rem-specific.approx}) to each of
$E$, $I$ and $J$ independently. Note that $g-a \geq d^{-\gb_1} g$
and $h-b \geq d^{-\gb_1} h$ follow from the assumptions on $\gd$
in Lemma \ref{lem-main} (recall $|[E]|=a\leq M/2$ and
$|[I]|=|I|\leq |E|\leq M/2$), justifying the first two
applications of Lemma \ref{lem-basic.approx} and the dependence of
$c_1$ on $\gb_1$ alone.

For the third application, note that if $b' \leq M/2$ we have
$h'-b' \geq d^{-\gb_1} h'$ (recall $|[J]|=|J|$). If $b' > M/2$,
then $|\cO\setminus N(J)| \leq M/2$. Since $[J]=J$ we also have
$[\cO\setminus N(J)]=\cO\setminus N(J)$ and $N(\cO\setminus
N(J))=\cE\setminus J$ and so (by the bound on $\gd$)
$(M-b')-(M-h') \geq d^{-\gb_1}(M-b')$. Using $h' \geq b'$ it
follows that $h'-b' \geq d^{-\gb_1}(M-h')$. But since $N(J) \cap
N(E) = \emptyset$ we have $h' \leq M-g \leq M-a \leq (1-\gb_2)M$
and so $h'-b' \geq d^{-\gb_1}(1/(1-\gb_2)-1)h' \geq d^{-c} h'$,
where the constant $c$ depends on both $\gb_1$ and $\gb_2$. \qed

\medskip

Before going on to the final step in the proof of Lemma
\ref{lem-main}, we need the following simple inequalities
(\cite[Lemma 3.1]{GalvinTetali} in the case $\psi=\sqrt{d}$). If
$(F,S,P,Q,P',Q')$ is an approximation for $(E,O) \in
\cH(a,g,b,h,b',h')$ then for suitably large $d$
\begin{equation} \label{inq-bounding.s.by.f}
|S| \leq |F| + \frac{3(g-a)}{\sqrt{d}},~~~~~|Q| \leq |P| +
\frac{3(h-b)}{\sqrt{d}}~~~~~\mbox{and}~~~~~ |Q'| \leq |P'| +
\frac{3(h'-b')}{\sqrt{d}}.
\end{equation}
Bearing this and the fact that $g-a \geq \gd g$ in mind, Lemma
\ref{lem-main} is implied by Lemma \ref{lem-adv.approx} and the
following reconstruction lemma.
\begin{lemma} \label{lem-reconstruction}
Let $G$, $a$, $g$, $b$, $h$, $b'$ and $h'$ be as in Lemma
\ref{lem-main}. There are constants $c_1=c_1(\gb_1)>0$ and
$c_2=c_2(\gb_1,\gb_2)>0$ such that for each $(F,S,P,Q,P',Q')
\subseteq \cO \times \cE \times \cE \times \cO \times \cO \times
\cE$ satisfying (\ref{inq-bounding.s.by.f}) there are at most
$$
\exp_2 \left\{M-b-b'-\frac{c_1(g-a)}{\log d}-\frac{c_1(h-b)}{\log
d}-\frac{c_2(h'-b')}{\log d}\right\}
$$
pairs $(E,O) \subseteq \cE \times \cO$ satisfying
\begin{equation} \label{containment.relations}
F \subseteq N(E),~~~S \supseteq [E],~~~P \subseteq N(I),~~~Q
\supseteq I,~~~P' \subseteq N(J)~~~\mbox{and}~~~Q' \supseteq J.
\end{equation}
\end{lemma}

\noindent {\em Proof: }For notational convenience, write $t$ for
$g-a$, $s$ for $h-b$ and $s'$ for $h'-b'$. Say that $S$ is {\em
tight} if $|S| < g - c_1't/\log d$ and {\em slack} otherwise, that
$Q$ is {\em tight} if $|Q| < b + c_1's/\log d$, and {\em slack}
otherwise, and that $Q'$ is {\em tight} if $|Q'| < b' + c_2's'/\log
d$, and {\em slack} otherwise, where $c_1'=c_1'(\gb_1)>0$ and
$c_2'=c_2'(\gb_1,\gb_2)>0$ are constants that will be specified
presently.

We now describe a procedure which, for input $(F,S,P,Q,P',Q')$
satisfying (\ref{inq-bounding.s.by.f}), produces an output $(E,O)$
which satisfies (\ref{containment.relations}). The procedure
involves a sequence of choices, the nature of the choices depending
on whether $S$, $Q$ and $Q'$ are tight or slack.

We begin by identifying a subset $D$ of $E$ which can be specified
relatively cheaply: if $Q$ is tight, we pick $I \subseteq Q$ with
$|I|=b$ and take $D=N(I)$; if $Q$ is slack, we simply take $D=P$
(recalling that $P \subseteq N(I) \subseteq E$).

If $S$ is tight, we complete the specification of $E$ by choosing $E
\setminus D \subseteq S \setminus D$. If $S$ is slack, we first
complete the specification of $N(E)$ by choosing $N(E) \setminus F
\subseteq N(S) \setminus F$. We then complete the specification of
$E$ by choosing $E \setminus D \subseteq [E] \setminus D$ (noting
that we do know $[E] \setminus D$ at this point).

Next we turn to the specification of $O$. As with $E$, we begin by
identifying a subset $D'$ of $O$: if $Q'$ is tight, we pick $J
\subseteq Q'$ with $|J|=b'$ and take $D'=N(J)$; if $Q'$ is slack,
we simply take $D'=P'$. From here, we complete the specification
of $O$ by choosing $O\setminus D' \subseteq \cO \setminus (N(E)
\cup D')$ (recall that $E \not \sim O$).

This procedure produces all pairs $(E,O)$ satisfying
(\ref{containment.relations}). Before bounding the number of
outputs, we gather together some useful observations.

First note that as established in the proof of Lemma
\ref{lem-adv.approx} we have
\begin{equation} \label{lower.bounds.on.g-a.etc}
g-a \geq d^{-\gb_1}g, ~~~~~ h-b \geq d^{-\gb_1}h ~~~~~ \mbox{and}
~~~~~ h'-b' \geq d^{-c}h'
\end{equation}
where the constant $c>0$ depends on both $\gb_1$ and $\gb_2$,
while from (\ref{inq-bounding.s.by.f}) we have
\begin{equation} \label{sqq'}
|S| \leq 2g, ~~~~~ |Q| \leq 2h ~~~~~ \mbox{and} ~~~~~ |Q'| \leq 2h'
\end{equation}
for suitably large $d$.

If $Q$ is tight then there are at most
\begin{equation} \label{choicesfordqsmall}
\sum_{i \leq c_1's/\log d}{|Q| \choose |Q|-i} \leq \sum_{i \leq
c_1's/\log d}{2h \choose i} \nonumber \leq 2^{s/2}
\end{equation}
possibilities for $D$ (for sufficiently small choice of the
constant $c_1'$, depending on $\gb_1$), and in this case $|D|=h$.
Here we are using (\ref{inq-binomial.2}) and
(\ref{lower.bounds.on.g-a.etc}). If $Q$ is slack there is just one
possibility for $D$, and in this case (using
(\ref{inq-bounding.s.by.f}))
\begin{equation} \label{qlarge}
|D| > b+c_1's/\log d  -3s/\sqrt{d} \geq b+c_1's/2\log d
\end{equation}
for suitably large $d$.

Similarly if $Q'$ is tight then there are at most $2^{s'/2}$
possibilities for $D'$ (for suitably small $c_2'$ depending on
both $\gb_1$ and $\gb_2$; here we use
(\ref{lower.bounds.on.g-a.etc})), and in this case $|D'|=h'$;
while if $Q'$ is slack there is just one possibility for $D'$, and
in this case $|D'|
> b'+c_2's'/2\log d$ for suitably large $d$.

If $S$ is slack then (\ref{inq-bounding.s.by.f}) implies $|N(E)
\setminus F| < 2c_1't/\log d$ and since $|N(S)\setminus F| \leq d|S|
\leq 2dg$ (see (\ref{sqq'})) the number of possibilities for $N(E)
\setminus F$ is at most
\begin{equation}
\sum_{i<2c_1't/\log d}{2gd \choose i} \leq 2^{t/2}
\label{choicesforglessfslarge}
\end{equation}
for suitable small $c_1'$ depending on $\gb_1$ (here again we use
(\ref{lower.bounds.on.g-a.etc})).

Now we assume that $d$, $c_1'$ and $c_2'$ are suitably chosen so
that all the previously made observations hold. We bound the number
of outputs of the procedure, considering first the four cases
determined by whether $S$ and $Q$ are slack or tight, and then
considering the two cases of whether $Q'$ is slack or tight. If $S$
and $Q$ are both tight then the number of possibilities for $E$ is
at most
\begin{equation} \label{ssmallqsmall}
\exp_2\{(s/2)+(g - c_1't/\log d-h)\}  \leq   \exp_2\{g - c_1't/\log
d-b-s/2\}.
\end{equation}
(The first term in the exponent on the left-hand side corresponds
to the choice of $D$ (using (\ref{choicesfordqsmall})), and the
second to the choice of $E\setminus D \subseteq S\setminus D$
(note that since $S$ and $Q$ are both tight, $|S\setminus D| \leq
g - c_1't/\log d-h$).

If $S$ is tight and $Q$ is slack then the total is at most
\begin{equation} \label{ssmallqlarge}
\exp_2\{g - c_1't/\log d-b-c_1's/2\log d\}.
\end{equation}
(Here there is no choice for $D$, and the exponent corresponds to
the choice of $E\setminus D \subseteq S\setminus D$ (using
(\ref{qlarge})).)

If $Q$ is tight then $|[E]\setminus D| = a-h$, so that if $S$ is
slack (and $Q$ tight) then the number of possibilities for $E$ is at
most
\begin{equation} \label{slargeqsmall}
\exp_2\{(s/2)+(t/2)+(a-h)\} \leq \exp_2\{g-t/2-b-s/2\}.
\end{equation}
(The first term in the exponent on the left-hand side corresponds
to the choice of $D$ (using (\ref{choicesfordqsmall})), the second
to the choice of $N(E)\setminus F$ (using
(\ref{choicesforglessfslarge})) and the third to the choice of
$E\setminus D$.)

If $Q$ is slack then $|[E]\setminus D| \leq a- b-c_1's/2\log d$ (see
(\ref{qlarge})), so that if $S$ and $Q$ are both slack the number of
possibilities for $E$ is at most
\begin{equation} \label{slargeqlarge}
\exp_2\{(t/2)+(a-b-c_1's/2\log d)\} \leq \exp_2\{g-t/2-b-c_1's/\log
d\}.
\end{equation}
(The first term in the exponent on the left-hand side corresponds to
the choice of $N(E)\setminus F$ and the second to the choice of
$E\setminus D$.)

Now we consider the number of choices for $O$, given our choice of
$E$. Note that $O \subseteq (\cO \setminus N(E))$, a set of size
$M-g$. If $Q'$ is tight then the number of possibilities for $O$ is
at most
\begin{equation} \label{q'small}
\exp_2\{(s'/2)+(M-g-h')\}  \leq   \exp_2\{M-g - b' -s'/2\}.
\end{equation}
(The first term in the exponent on the left-hand side corresponds
to the choice of $D'$, and the second to the choice of $O\setminus
D' \subseteq \cO\setminus (N(E) \cup D')$.)

Finally if $Q'$ is slack then the number of possibilities for $O$ is
at most
\begin{equation} \label{q'large}
\exp_2\{M-g-b'-c_2's'/2\log d\}.
\end{equation}
(Here there is no choice for $D'$, and the exponent corresponds to
the choice of $O\setminus D' \subseteq \cO\setminus (N(E) \cup
D')$.)

Combining (\ref{ssmallqsmall}), (\ref{ssmallqlarge}),
(\ref{slargeqsmall}), (\ref{slargeqlarge}), (\ref{q'small}) and
(\ref{q'large}) we obtain the lemma. \qed


\begin{thebibliography}{99}

\bibitem{AldousFill}
D. Aldous and J. Fill, {\em Reversible Markov Chains and Random
Walks on Graphs}, monograph in preparation, available at
\url{http://stat-www.berkeley.edu/users/aldous/RWG/book.html}.

\bibitem{BenjaminiHaggstromMossel}
I. Benjamini, O. H\"aggstr\"om and E. Mossel, On random graph
homomorphisms into ${\mathbb Z}$, {\em J. Combinatorial Th. (B)}
{\bf 78} no. 1 (2000), 86--114.

\bibitem{Bollobas2}
B. Bollob\'as, {\em Extremal Graph Theory}, Academic Press, New
York, 1978.

\bibitem{Bollobas}
B. Bollob\'as, {\em Modern Graph Theory}, Springer, New York, 1998.

\bibitem{Bollobas3} B. Bollob\'as, {\em Random Graphs}, Cambridge University Press,
Cambridge, 2001.

\bibitem{BorgsChayesFriezeKimTetaliVigodaVu}
C. Borgs, J. Chayes, A. Frieze, J. Kim, P. Tetali, E. Vigoda, V. Vu,
Torpid Mixing of some Monte Carlo Markov Chain algorithms in
Statistical Physics, {\em Proc. IEEE FOCS '99}, 218--229.

\bibitem{BubleyDyer}
R. Bubley and M. Dyer, Path coupling: a technique for proving rapid
mixing in Markov chains, {\em Proc. IEEE FOCS '97}, 223--231.

\bibitem{Chernoff} H. Chernoff, A measure of asymptotic efficiency for tests of
a hypothesis based on the sum of observations, {\em Ann. Math.
Statistics} {\bf 23} (1952), 493--507.

\bibitem{Diestel}
R. Diestel, {\em Graph Theory}, Springer, New York, 2005.

\bibitem{DyerFriezeJerrum}
M. Dyer, A. Frieze and M. Jerrum, On counting independent sets in
sparse graphs, {\em SIAM J. Comp.} {\bf 31} (2002), 1527--1541.

\bibitem{DyerGreenhillMolloy}
M. Dyer, C. Greenhill, M. Molloy, Very rapid mixing of the Glauber
dynamics for proper colorings on bounded degree graphs, {\em Random
Struc. \& Alg.} {\bf 20} (2002), 98--114.

\bibitem{Galvin}
D. Galvin, On homomorphisms from the Hamming cube to $\Z$, {\em
Isr. J. Math.} {\bf 138} (2003), 189--213.

\bibitem{GalvinRandall}
D. Galvin and D. Randall, Torpid Mixing of Local Markov Chains on
$3$-Colorings of the Discrete Torus, {\em Proc. ACM--SIAM SODA '07},
376--384.

\bibitem{GalvinTetali}
D. Galvin and P. Tetali, Slow mixing of Glauber dynamics for the
hard-core model on the Hamming cube, {\em Random Structures \& Alg.}
{\bf 28} (2006) 427-443.

\bibitem{Jerrum}
M. Jerrum, A very simple algorithm for estimating the number of
$k$-colourings of a low-degree graph, {\em Random Struc. \& Alg.}
{\bf 7} (1995), 157--165.

\bibitem{JerrumSinclair}
M. Jerrum and A. Sinclair, Conductance and the rapid mixing property
for Markov chains: the approximation of the permanent resolved, {\em
Proc. ACM STOC '88}, 235--243.

\bibitem{Kahn}
J. Kahn, Range of cube-indexed random walk, {\em Isreal J. Math.}
{\bf 124} (2001) 189--201.

\bibitem{KorshunovSapozhenko}
A. Korshunov and A. Sapozhenko, The number of binary codes with
distance $2$, {\em Problemy Kibernet.} {\bf 40} (1983), 111--130.
(Russian)

\bibitem{LuczakVigoda}
T. \L uczak and E. Vigoda, Torpid mixing of the
Wang-Swendsen-Koteck\'y algorithm for sampling colorings, {\em J.
Discrete Alg.} {\bf 3} no. 1 (2005), 92--100.

\bibitem{Molloy}
M. Molloy, Very rapidly mixing Markov chains for $2\gD$-colourings
and for independent sets in a $4$-regular graph, {\em Random Struc.
\& Alg.} {\bf 18} (2001), 101--115.

\bibitem{MontenegroTetali}
R. Montenegro and P. Tetali, Mathematical aspects of mixing times in
Markov chains, Foundations and Trends in Theoretical Computer
Science {\bf 1} no. 3 (2006), 237--354.

\bibitem{Randall}
D. Randall, Mixing, {\em Proc. IEEE FOCS '03}, 4--15.

\bibitem{Randall2}
D. Randall, personal communication.

\bibitem{SalasSokal}
J. Salas and A. Sokal, Absence of phase transition for
antiferromagnetic Potts models via the Dobrushin uniqueness theorem,
{\em J. Stat. Phys.} {\bf 86} (1997), 551--579.

\bibitem{Sapozhenko}
A. Sapozhenko, On the number of connected subsets with given
cardinality of the boundary in bipartite graphs, {\em Metody
Diskret. Analiz.} {\bf 45} (1987), 42--70.  (Russian)

\bibitem{Sapozhenko2}
A. A. Sapozhenko, The number of antichains in ranked partially
ordered sets, {\em Diskret. Mat.} {\bf 1} (1989), 74--93. (Russian;
translation in Discrete Math. Appl. 1 no. 1 (1991), 35--58)

\bibitem{Sokal1}
A. Sokal, Chromatic Polynomials, Potts Models and All That, {\em
Physica} {\bf A}279 (2000), 324--332.

\bibitem{Sokal2}
A. Sokal, A Personal List of Unsolved Problems Concerning Lattice
Gases and Antiferromagnetic Potts Models, {\em Markov Process.
Related Fields} {\bf 7} (2001), 21--38.

\bibitem{Vigoda}
E. Vigoda, Improved bounds for sampling colorings, {\em J. Math.
Phys.} {\bf 41} (2000), 1555--1569.

\end{thebibliography}
\end{document}